\documentclass[12pt,a4paper]{article}

\usepackage{latexsym,amsmath,amsthm,amssymb,mathrsfs}
\usepackage[all]{xy}
\usepackage{textcomp}
\usepackage{graphicx}

\begin{document}

\author{Carlotta Gazzola, Ezio Venturino\\
Dipartimento di Matematica ``Giuseppe Peano'',\\
Universit\`a di Torino, \\
via Carlo Alberto 10, 10123 Torino, Italy\\
E-mail: ezio.venturino@unito.it, 
}
\date{}
\title{Dynamics of different topological configurations in a three-patch metapopulation model}
\maketitle

\begin{abstract}
The possible dynamics of an ecosystem with three interconnected patches among which one population
can migrate are investigated, paying attention to the consequences that possible disruption in the
communicating paths can cause.
\end{abstract}


AMS codes 92D30, 92D25, 92D40


\section{Introduction}

Metapopulations are a tool for the theoretical investigation of fragmented landscapes, \cite{HG,HT}. Whether the fragmentation is due
to human activity or natural causes like landslides, fires and so on, it might represent for the animal populations living
in the interested ecosystem a possible danger. Basically, this framework consists in formulating models in which to the
local population dynamics interpatch migrations are suitably added. Mathematical models are needed even by field ecologists,
in view of the difficulty of gathering migration data.
Persistence of wild populations in heterogeneous environments is of particular concern for conservationist biologists,
\cite{W96,W97}.

Metapopulation dynamics has been successfully applied to model 
for instance the mountain sheep ({\it{Ovis canadensis}}), \cite{FB},
or the spotted owl 
({\it {Strix occidentalis}}), \cite{GH}.

The study of communities living in separate environments that are connected by possible migration paths
has led to some counterintuitive results, such as the fact that
the global population thrives, while locally in some patches becoming extinct,
\cite{C03,GG,H85,Wu}. In this setting, also recent works on the role of niches as safety refuges can be
accounted for, \cite{Collings1995,Gonzalez2003,Gonzalez2012,Gonzalez2012Prey,Ruxton1995short}.

More recently, \cite{EV_11}, in the framework of modelling heterogenous habitats, also the influence of diseases have been introduced.
The effects of diseases on communities are a fact in nature and therefore also demographic models cannot ignore them. In fact, since
a quarter of a century now, models accounting for interacting populations where also epidemics spread have been proposed and investigated, see Chapter 7 of \cite{MPV} and the papers
\cite{HF,BC,EV94,EV95,CA,AAMC,EV01,EV02,EV:07}.

In general, the mathematical models that are introduced in this context can contain several patches and are usually
analysed for the existence of the equilibria and possibly for their stability. In this paper we want to consider a rather
simple system, composed of three patches that are joined together by connecting directed paths. Our aim is to investigate
how its behavior changes once some of these connections are broken, whether accidentally or, as mentioned above, due to human
artifacts that partly or entirely disrupt these communications between these habitats. 

The paper is organized as follows. The general model with all possible connections between the three environments is presented
and analysed in the next Section, finding its possible equilibria and studying their stability.
Section \ref{sec:br} contains the models in which some of the paths become unavailable for the animals interpatch movement.
A final discussion of the ecological implications concludes the paper.

\section{The general model}

We consider here an environment made out of 3 patches that are interconnected for migrations of a population $P$, as depicted in the following diagram.
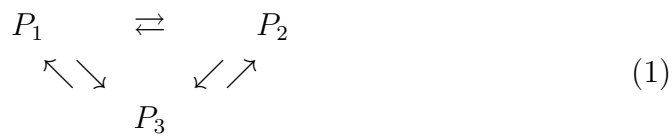
\begin{eqnarray}
P_1 \qquad & \rightleftarrows & \qquad P_2\nonumber\\ \label{fig_sist4}
\nwarrow \searrow &&  \swarrow \nearrow \\
&P_3&\nonumber
\end{eqnarray}
The size of each subpopulation in each patch $k$ is denoted by $P_k$,
$k=1,2,3$. Assuming that migrations from each patch are possible in all directions toward both other patches, the model is
\begin{eqnarray}\label{sist4}
\dot{P_1}=r_1P_1\left(1-\displaystyle\frac{P_1}{k_1}\right)+m_{12}P_2 +m_{13}P_3-m_{21}P_1-m_{31}P_1\\ \nonumber
\dot{P_2}=r_2P_2\left(1-\displaystyle\frac{P_2}{k_2}\right)-m_{32}P_2 -m_{12}P_2+m_{21}P_1+m_{23}P_3\\ \nonumber
\dot{P_3}=r_3P_3\left(1-\displaystyle\frac{P_3}{k_3}\right)+m_{32}P_2-m_{13}P_3-m_{23}P_3+m_{31}P_1
\end{eqnarray}
All the parameters are implicitly assumed to be nonnegative.
Each equation describes the population dynamics in each patch. Each subpopulation reproduces logistically,
with parameters that are environment-dependent, namely net reproduction rate
$r_i$ and carrying capacity $k_i$, $i=1,2,3$. In addition, denoting by $m_{ij}$ the migration rates from
patch $j$ into patch $i$, we assume that movements in between different patches depend directly from the
population level in the outgoing patch.

For later stability analysis purposes, it is convenient also to consider the Jacobian $J$ of (\ref{sist4}),
\begin{equation}\label{J_4}
\displaystyle
\left[
\begin{array}{ccc}
J_{11}& m_{12} & m_{13}\\ 
m_{21}& J_{22} & m_{23} \\
m_{31}& m_{32}& J_{33}
\end{array}
\right]
\end{equation}
with
\begin{eqnarray*}
J_{11}=r_1-2\frac{r_1}{k_1}P_1-m_{21}-m_{31},\\
J_{22}=r_2-2\frac{r_2}{k_2}P_2-m_{12}-m_{32},\\
J_{33}=r_3-2\frac{r_3}{k_3}P_3-m_{13}-m_{23}.
\end{eqnarray*}

\subsection{Feasible equilibria}
There are only two possible equilibria, the origin at which the ecosystem disappears, and possibly the coexistence equilibrium,
in which all the patches are populated. We now prove its existence.

Solve for $P_3$ the first two equations of (\ref{sist4}) obtaining two surfaces:
\begin{eqnarray*}
P^{(1)}_3=\frac{-1}{m_{13}}\left[r_1P_1\left(1-\frac{P_1}{k_1}\right)+m_{12}P_2 -m_{21}P_1-m_{31}P_1\right]\\
P^{(2)}_3=\frac{-1}{m_{23}}\left[r_2P_2\left(1-\frac{P_2}{k_2}\right)-m_{32}P_2 -m_{12}P_2+m_{21}P_1\right]
\end{eqnarray*}
The restriction of $P^{(1)}_3$ to the $P_2-P_3$ coordinate plane is a straight line with negative slope through the origin.
Instead, the restriction of $P^{(1)}_3$ to the $P_1-P_3$ coordinate plane is a convex parabola
through the origin.
This parabola has another zero at the point $P_1^{(0,1)}=k_1r_1^{-1}[r_1-(m_{21}+m_{31})]$. Hence a feasible branch
emanates from this latter point, when $P_1^{(0,1)}\ge 0$,
or from the origin in the opposite case.
A similar result holds for $P^{(2)}_3$, where in this case the root is
$P_1^{(0,2)}=k_2r_1^{-2}[r_2-(m_{12}+m_{32})]$
Note that $P_1^{(0,1)}\ge 0$ and $P_1^{(0,2)}\ge 0$ hold when respectively the following conditions
are satisfied
\begin{equation}\label{zeros_restr}
r_1 \ge m_{21}+m_{31}:= M_1, \quad  r_2\ge m_{12}+m_{32}:= M_2.
\end{equation}

To better study the problem, however, we consider the intersections of these surfaces with the horizontal
planes, $P_3=h\ge 0$. Again, two parabolae are found,
\begin{eqnarray*}
\sigma_h : P_2=\frac{-1}{m_{12}}\left[r_1P_1\left(1-\frac{P_1}{k_1}\right)-\left( m_{21}+m_{31}\right) P_1 + m_{13}h\right]\\
\rho_h : P_1=\frac{-1}{m_{21}}\left[r_2P_2\left(1-\displaystyle\frac{P_2}{k_2}\right)-\left( m_{32}+m_{12}\right) P_2+m_{23}h\right]
\end{eqnarray*}
with $\sigma_h$ being a convex function of $P_2$ and $\rho_h$ a convex function of $P_1$.
Both have nonpositive values at the origin, so that their two roots are of opposite signs.
%
Hence a feasible branch
emanates from the positive root, when $h> 0$,
or from the origin where this root degenerates for $h=0$.

Both parabolae have only one branch that lies in the feasible orthant. Hence, the two
curves must meet at exactly one point $Q_h$, with nonnegative coordinates. In particular on the $P_3=0$ plane, i.e. for $h=0$, 
$Q_0$ at worst could coincide with the origin in case (\ref{zeros_restr}) both do not hold.
Since $h$ is arbitrary, it follows that the two surfaces
$P^{(1)}_3$ and $P^{(2)}_3$ meet along a line $\ell$ in the positive orthant:
$$
\ell= \left\{ Q_h\equiv \sigma_h \cap \rho_h : \forall h\ge 0 \right\}.
$$

We now consider the surfaces $\Sigma_{\pm}:P_3^{\pm}\equiv P_3^{\pm}(P_1,P_2)$ originating from the third equation of (\ref{sist4}),
given explicitly by
$$
P_3^{\pm}=\frac {k_3}{2r_3}\left[ r_3-m_{13}-m_{23}\pm\sqrt{(r_3-m_{13}-m_{23})^2+4\frac{r_3}{k_3}(m_{32}P_2+m_{31}P_1)}\right].
$$
Since the term under the square root exceeds the one outside it, 
for every possible value of $P_1\ge0$ and $P_2\ge0$ the surface $\Sigma_+$ is always nonnegative,
while $\Sigma_-$ is always nonpositive.
Hence the intersection of $\ell$ and $\Sigma_+$ in the first orthant is always guaranteed, which
provides the unique feasible coexistence equilibrium.

In summary, we have shown the following result.

\vspace{0.3cm}
{\textbf{Theorem} The coexistence equilibrium always exists.
}

\subsection{Equilibria stability}

The stability of the equilibria can be assessed rather simply by using Descartes' rule of signs
on the characteristic equation. The latter is the cubic
$$
\lambda ^3 - {\textbf{tr}}(J) \lambda ^2 + M_J \lambda - \det (J)=0,
$$
where $M_J$ represents the sum of the principal minors of $J$ of order 2. To have all negative roots,
we need the conditions
$$
{\textbf{tr}}(J)<0, \quad M_J>0, \quad \det (J)<0.
$$
Let
$$
\Pi_i = r_i \left( 1- \frac 2{k_i} P_i\right)
$$
Explicitly, they become
\begin{equation}\label{traceJ}
\sum _{i=1}^3 \Pi_i< m_{21}+m_{31}+m_{12}+m_{32}+m_{13}+m_{23}
\end{equation}
for the trace,
\begin{eqnarray}\label{MJ}
\Pi_1\Pi_2+\Pi_1\Pi_3+\Pi_2\Pi_3+m_{31}m_{32}+m_{21}m_{32}+m_{13}m_{31}\\ \nonumber
+m_{12}m_{23}+m_{31}m_{23}+m_{13}m_{21}+
m_{21}m_{13}+m_{12}m_{23}+m_{32}m_{13}\\ \nonumber
>\Pi_1(m_{32}+m_{12}+m_{13}+m_{23})+
\Pi_2(m_{31}+m_{21}+m_{13}+m_{23})\\ \nonumber
+\Pi_3(m_{21}+m_{31}+m_{12}+m_{32})
\end{eqnarray}
for $M_J$ and finally for the determinant we have
\begin{eqnarray}\label{detJ}
&&\Pi_1\Pi_2\Pi_3 + \Pi_1(m_{12}m_{13}+m_{12}m_{23}+m_{32}m_{13})\\ \nonumber
&&+\Pi_2(m_{21}m_{31}+m_{21}m_{23}+m_{31}m_{23})
+\Pi_3(m_{21}m_{32}+m_{31}m_{12}+m_{31}m_{32})\\ \nonumber
&&> \Pi_1\Pi_2 (m_{13}+m_{23}) 
+\Pi_1\Pi_3 (m_{12}+m_{32})
+\Pi_2\Pi_3 (m_{21}+m_{31}).
\end{eqnarray}
For the origin, note the simplification $\Pi_i=r_i$.
For the coexistence equilibrium, the above stability conditions are more involved to assess.
Numerical simulations however reveal its stability. 


\section{The models with some broken paths}\label{sec:br}

There are several situations that can arise, when due to human artifacts or some natural catastrophic events some
of the connecting paths become unavailable for the population migrations.
We avoid to consider the situations in which one or all the
patches become isolated. In this case indeed the isolated subpopulation would thrive independently of the others, due to the intrinsic
resources represented by each logistic model in the formulation of (\ref{sist4}), and the
remaining configuration is
simply given by two possibly connected patches, and therefore it is very easy to analyse.
There are thus nine possible situations.

In fact, we can remove one (directed) path between any of the 3 patches in just one way. Combinatorically
indeed it does not make any difference among which nodes we decide to break the connection and furthermore also the direction of the
removed arc is immaterial, since by relabeling the nodes we would end up with the same situation. The result is represented 
in the picture below.

\begin{eqnarray}\label{2}
P_1 \quad \quad &\leftrightarrows& \quad P_2 \qquad \qquad (EX 2) = (5) \nonumber\\
\searrow \nwarrow&& \nearrow\\
&P_3& \nonumber
\end{eqnarray}

We can then remove two edges in several ways. From the same nodes, we get the following configuration

\begin{eqnarray}\label{0}
P_2 \qquad & 
& \qquad P_3\nonumber\\ 
\nwarrow \searrow &&  \swarrow \nearrow \\
&P_1&\nonumber
\end{eqnarray}

If they are removed from different connected nodes, there are three alternatives: either the
node that is connected with both the other nodes by just one arc
has the two edges one outgoing and one incoming,
or both outgoing, or both incoming. The pictures below will better illustrate these 3 situations.

\begin{eqnarray}\label{3}
P_3 \quad &\leftrightarrows& \quad P_2 \qquad \qquad (EX 3) \nonumber\\
\searrow && \nearrow \\
&P_1&\nonumber
\end{eqnarray}

\begin{eqnarray}\label{7}
P_3 \quad \quad &\leftarrow& \quad P_2  \qquad \qquad (EX 7) \nonumber\\
\searrow \nwarrow && \swarrow \\
&P_1&\nonumber
\end{eqnarray}

\begin{eqnarray}\label{8}
P_3 \quad &\leftrightarrows& \quad P_2 \qquad \qquad (EX 8) \nonumber\\
\searrow && \swarrow \\
&P_1&\nonumber
\end{eqnarray}

Next, we can remove 3 edges.
If we remove one edge from each pair of nodes, the
only alternative is the way in which the orientation is considered. We can either remove all the edges in the same direction,
but in such case which direction is immaterial, by a suitable relabeling of the nodes,
or one edge is removed in one direction and the remaining two in the opposite one; again due to symmetries this leads to
just one configuration. These alternatives are depicted below.

\begin{eqnarray}\label{1}
P_3 \quad &\leftarrow& \quad P_2 \qquad \qquad (EX 1) \nonumber\\
\searrow && \nearrow \\
&P_1&\nonumber
\end{eqnarray}

\begin{eqnarray}\label{6}
P_3 \quad &\leftarrow& \quad P_2 \qquad \qquad (EX 6) \nonumber\\
\searrow && \swarrow \\
&P_1&\nonumber
\end{eqnarray}
If we remove 2 edges connecting the same nodes and remove another one, apart from symmetries there are only two configurations
possible, namely

\begin{eqnarray}\label{2bis}
P_1 \quad \quad & & \quad P_2 \qquad \qquad (EX 2\ NEW) = (5) \nonumber\\
\searrow \nwarrow&& \nearrow\\
&P_3& \nonumber
\end{eqnarray}

\begin{eqnarray}\label{7bis}
P_1 \quad \quad & & \quad P_2  \qquad \qquad (EX 7\ NEW) \nonumber\\
\searrow \nwarrow && \swarrow \\
&P_3&\nonumber
\end{eqnarray}

Finally, we can remove 4 edges and no more, otherwise at least one node will be disconnected from the other ones.
The possible system configurations are as follows:
\begin{eqnarray}\label{n1}
P_1 \quad \rightarrow \quad &P_2 \quad \rightarrow \quad P_3
\end{eqnarray}
\begin{eqnarray}\label{n2}
P_1 \quad \rightarrow \quad P_2 \quad \leftarrow \quad P_3
\end{eqnarray}
\begin{eqnarray}\label{n3}
P_1 \quad \leftarrow \quad P_2 \quad \rightarrow \quad P_3
\end{eqnarray}

For the analysis of all these models with broken paths, we must set some of the migration rates to zero.
We will then investigate whether new equilibria arise, and, if possible, whether the stability of the
origin and of coexistence are altered in the new configurations.

\subsection{The models (\ref{2}), (\ref{0}), (\ref{3}), (\ref{1})}
In general, in these three models, the only
equilibria are those of the general model, as no other ones can arise.
For (\ref{0}) no changes are necessary in the proof of the coexistence equilibrium.
In the other cases, the proof however requires some attention.

Specifically, consider first (\ref{2}) where $m_{23}=0$. Note that in this case the approach used to show the
existence of the equilibrium fails. We can instead solve the first and third equilibrium equations
for $P_2$ and intersect the corresponding surfaces $P_2^{(1)}(P_1,P_3)$ and $P_2^{(2)}(P_1,P_3)$ with
the planes $P_2=h$ to get
\begin{eqnarray*}
\alpha_h : P_3=\frac{-1}{m_{13}}\left[r_1P_1\left(1-\frac{P_1}{k_1}\right)-\left( m_{21}+m_{31}\right) P_1 + m_{12}h\right]\\
\beta_h : P_1=\frac{-1}{m_{31}}\left[r_3P_3\left(1-\displaystyle\frac{P_3}{k_3}\right)-m_{13} P_3+m_{32}h\right].
\end{eqnarray*}
These parabolae are seen to intersect with each other along a line that itself intersects the remaining
surface originating from the second equilibrium equation, as done in the proof of the Theorem.

For (\ref{3}) instead we have
$m_{31}=m_{12}=0$, solving the first and second equilibrium equation as done in the Theorem,
for $P_3^{(1)}$ we obtain a parabolic cyclinder. Its intersection with $P_3=h$ gives a straight line,
and the latter always meets the parabola obtained intersecting $P_3^{(2)}$ with $P_3=h$, as in the
proof of the Theorem, and existence follows accordingly.

Finally for (\ref{1}) in which $m_{31}=m_{12}=m_{23}=0$,
solving all equilibrium equations
we obtain always parabolic cylinders, with each axis parallel to a different coordinate axis:
\begin{eqnarray*}
\gamma_1 : P_3=\frac{-1}{m_{13}}\left[r_1P_1\left(1-\frac{P_1}{k_1}\right)-m_{21} P_1\right],\\
\gamma_2 : P_1=\frac{-1}{m_{21}}\left[r_2P_2\left(1-\displaystyle\frac{P_2}{k_2}\right)-m_{32}P_2\right],\\
\gamma_3 : P_2=\frac{-1}{m_{32}}\left[r_3P_3\left(1-\displaystyle\frac{P_3}{k_3}\right)-m_{13} P_3\right].
\end{eqnarray*}
Intersecting $\gamma_1$ with the plane $P_3=h$ gives a straight line $L_1$ parallel to the $P_2$ axis
and $L_1\cap \gamma_2$ gives a point, and therefore we get the line in space
parametrized by $Q_{12}^h(k^{h},\ell^{h},h)$, which itself, as in the Theorem, must intersect the $\gamma_3$ surface,
thereby providing the coexistence equilibrium.

Furthermore, the stability analysis hinges
on the Jacobian (\ref{J_4}), in which the above conditions make some simplifications.
But even for the origin the stability conditions nevertheless remain quite involved.
More specifically, for all these models we find that (\ref{traceJ}) becomes sharper,
as the right hand side will become smaller.
The left hand side contains only demographic parameters, and explicitly no migration
rates. On the other hand it also contains the population levels at equilibrium,
which in turn depend on the migration rates. Changes in the latter could in
principle bring the equilibrium populations up or down and therefore influence stability as well.
Similar considerations hold for (\ref{MJ}) and (\ref{detJ}).

\vspace{0.3cm}
{\textbf{Remark}}. It is thus hard to state whether the stability conditions
will be easier or more difficult to be satisfied. These general considerations hold also for all
the other models with some broken paths, unless we explicitly present some further remarks.
\vspace{0.3cm}

\subsection{The model (\ref{2bis})}
Here we have $m_{12}=m_{21}=m_{32}=0$. These simplifications do not harm the proof of the Theorem,
so that the coexistence equilibrium is guaranteed to be feasible. In addition, however, they show
that the origin is certainly unstable, as one eigenvalue for this equilibrium is $J_{22}=r_2>0$.

In addition, this model
allows also the equilibrium $X=(0,k_2,0)$, for which one eigenvalue is explicit, $-r_2<0$,
and the remaining ones provide the stability conditions
\begin{eqnarray}\label{X_stab_mod7bis}
m_{31}+m_{32}+m_{13}>r_3+r_1, \\ \nonumber
(m_{31}-r_1)(m_{32}+m_{13}-r_3)>m_{13}m_{31}.
\end{eqnarray}

\subsection{The models (\ref{7}), (\ref{7bis})}

For (\ref{7}) we need to take $m_{21}=m_{23}=0$. It follows immediately that
$$
P_2^*=\frac {k_2}{r_2}\left( r_2-m_{12}-m_{32}\right).
$$
From the remaining equilibrium equations we discover that the two parabolae
\begin{eqnarray*}
P_3=\frac 1{m_{13}}\left[ m_{31}P_1-r_1P_1\left( 1-\frac {P_1}{k_1}\right)-m_{12} P_2^*\right],\\
P_1=\frac 1{m_{31}}\left[ m_{13}P_3-r_3P_3\left( 1-\frac {P_3}{k_3}\right)-m_{32} P_2^*\right],
\end{eqnarray*}
are convex, with negative value at 0, so that they
always meet in the first quadrant, thus providing the remaining components of the coexistence equilibrium.

Stability at coexistence comes just from the negativity of one explicit eigenvalue, providing
\begin{eqnarray}\label{stab_1_mod7bis}
m_{12}+m_{32}+2 \frac {r_2}{k_2}P_2^*>r_2,
\end{eqnarray}
while the Routh-Hurwitz conditions on the remaining minor hold always true,
\begin{eqnarray*}
\frac {r_1}{k_1}P_1^*+m_{12} \frac {P_2^*}{P_1^*}+m_{13}\frac {P_3^*}{P_1^*}
+\frac {r_3}{k_3}P_3^*+m_{32} \frac {P_2^*}{P_3^*}+m_{31}\frac {P_1^*}{P_3^*}>0,\\
\left( \frac {r_1}{k_1}P_1^*+m_{12} \frac {P_2^*}{P_1^*} \right)
\left( \frac {r_3}{k_3}P_3^*+m_{32} \frac {P_2^*}{P_3^*} \right)
+m_{13}\frac {P_3^*}{P_1^*} 
\left( \frac {r_3}{k_3}P_3^*+m_{32} \frac {P_2^*}{P_3^*} \right)\\
+m_{31}\frac {P_1^*}{P_3^*} 
\left( \frac {r_1}{k_1}P_1^*+m_{12} \frac {P_2^*}{P_1^*}\right)>0.
\end{eqnarray*}
Stability at the origin is obtained by
\begin{eqnarray}\label{stab_orig_mod7bis}
m_{12}+m_{32}>r_2, \quad
m_{13}+m_{31}>r_1+r_3, \quad
r_1r_3>r_1m_{13}+r_3m_{31}.
\end{eqnarray}

Further, for the model (\ref{7bis}) it is enough to set $m_{12}=0$ and all the above considerations still carry on
to this case.

In these models there is one more feasible equilibrium in addition to origin and coexistence, namely
$Q_1=(P_1^{Q},0,P_3^{Q})$.
To find it in both cases (\ref{7}) and (\ref{7bis}),
we solve the first and third equilibrium equations of (\ref{sist4}) to find
\begin{eqnarray}\label{parab_7}
P_3=\frac {P_1}{m_{13}} \left[ m_{31}-r_1 \left( 1-\frac {P_1}{k_1}\right) \right],\\
P_1=\frac {P_3}{m_{31}} \left[ m_{13}-r_3 \left( 1-\frac {P_3}{k_3}\right) \right].
\end{eqnarray}
These are two convex parabolae, with roots at the origin and respectively at the points
$P_1^{(a)}=k_1r_1^{-1}[r_1-m_{31}]$ and
$P_3^{(a)}=k_3r_3^{-1}[r_3-m_{13}]$. These points are nonnegative if the conditions
$r_1\ge m_{31}$ and $r_3\ge m_{13}$ hold.
When at least one of these conditions holds sharply, then an intersection between the parabolae is guaranteed in view of their convexity.
If instead both are equalities, then the parabolae are both tangent to the axes at the origin, and therefore
one intersection is the origin itself and another one exists also in this case.
When instead both are not satisfied, we need to compare the parabolae slopes at the origin to determine whether
an intersection between their feasible branches exists. We find
\begin{eqnarray*}
P'_3(P_1)=\frac 1{m_{13}} \left( m_{31}-r_1+2\frac {r_1}{k_1}P_1\right), \\
P'_1(P_3)=\frac 1{m_{31}} \left( m_{13}-r_3+2\frac {r_3}{k_3}P_3\right),
\end{eqnarray*}
and we must impose that
$$
P'_3(0)>(P^{-1}_1)'(0),
$$
in order to ensure that the branches meet in the first quadrant. This amounts to requiring
$$
r_1 (r_3 - m_{13}) - r_3 m_{31}>0,
$$
which is impossible, in view of the restrictions holding
in this situation
\begin{equation}\label{zeros_restr_2}
r_1<m_{31}, \quad r_3<m_{13}.
\end{equation}
Thus in this situation the point $Q_1$ is infeasible.

$Q_1$ is stable if the eigenvalue that is immediately found is negative,
entailing
\begin{equation}\label{Q1_stab_mod7bis}
r_2<m_{12}+m_{32},
\end{equation}
and for the model (\ref{7bis}) this condition simplifies since $m_{12}=0$. In fact,
the remaining Routh-Hurwitz conditions stemming from a 2 by 2 reduced Jacobian ${\widetilde J}$
are satisfied,
\begin{eqnarray*}
-{\mathrm{tr}}{\widetilde J}=\frac {m_{13}}{P_1^{Q}}P_3^{Q}+\frac {r_1}{k_1}P_1^{Q}
+\frac {m_{31}}{P_3^{Q}}P_1^{Q}+\frac {r_3}{k_3}P_3^{Q}>0, \\
\det {\widetilde J} =\frac {m_{13}r_3}{k_3P_1^{Q}}(P_3^{Q})^2+
\frac {m_{31}r_1}{k_1P_3^{Q}}(P_1^{Q})^2+
\frac {r_1r_3}{k_1k_3}P_1^{Q}P_3^{Q}>0
\end{eqnarray*}

Thus these models admit the origin, the patch-2-population-free point and coexistence as possible equilibria.

\subsection{The model (\ref{8})}

Here we set $m_{21}=m_{31}=0$. For coexistence the approach of the general case still works, it only
simplifies giving for $P_3^{(2)}$ a cylinder with axis parallel to the $P_1$ coordinate axis.
Stability at coexistence is guaranteed by the explicit eigenvalue,
\begin{equation}\label{Stab_8_coex}
k_1< 2P_1^*.
\end{equation}

{\textbf{Remark}}.
Note that this is an eigenvalue also for the origin, thereby providing its instability, in view of $k_1<0$.

The remaining Routh-Hurwitz conditions are always satisfied, since they reduce to
\begin{eqnarray*}
\frac {m_{23}}{P_2^*}P_3^*+\frac {r_2}{k_2}P_2^*
+\frac {m_{32}}{P_3^*}P_2^*+\frac {r_3}{k_3}P_3^*>0, \\
\frac {m_{23}r_3}{k_3P_3^*}(P_3^*)^2+
\frac {m_{32}r_2}{k_2P_3}(P_2^*)^2+
\frac {r_2r_3}{k_2k_3}P_2P_3>0.
\end{eqnarray*}

In addition to origin and coexistence, here we find also the equilibrium $M_2=(k_1,0,0)$
which is clearly always feasible.
For its stability, the Jacobian simplifies even further.
The first eigenvalue is $-r_1<0$, the Routh-Hurwitz 
conditions on the remaining ones give the stability conditions
\begin{eqnarray}\label{stab_82}
r_2+r_3<m_{12}+m_{32}+m_{13}+m_{23}, \\ \nonumber
(r_2-m_{12})(r_3-m_{13}) > (r_2-m_{12}) m_{23} + (r_3-m_{13}) m_{32}.
\end{eqnarray}
The former condition when becomes an equality gives rise to a Hopf bifurcation.
Explicitly, this occurs for
\begin{eqnarray}\label{hopf82}
r_2 =r_2^{\ddagger}= m_{13}+m_{23}+m_{32}+m_{12}-r_3.
\end{eqnarray}

\subsection{The model (\ref{6})}

Set $m_{21}=m_{31}=m_{23}=0$ in (\ref{sist4}).
Note that the origin in this model is unstable, since one eigenvalue is $J_{11}=r_1>0$.

Coexistence can be calculated explicitly, to give
\begin{eqnarray*}
P_1^*= \frac {k_1}2\left[ 1 + \sqrt{ 1 + \frac 4{k_1r_1}\left( m_{12} P_2^* + m_{13} P_3^* \right)} \right], \quad
P_2^*=\frac {k_2}{r_2}\left( r_2-m_{32}-m_{12}\right) , \\
P_3^*=\frac {k_3}{2r_3}\left[ r_3-m_{13}+ \sqrt{ \left( r_3-m_{13}\right)^2 +
\frac 4{k_3} r_3m_{32} P_2^*} \right] .
\end{eqnarray*}
Feasibility of the coexistence equilibrium is ensured just by
\begin{equation}\label{ce4}
r_2>m_{12}+m_{32}.
\end{equation}
The eigenvalues become
$$
J_{11}=-\frac 1{P^*_1} \left( m_{12}P^*_2 + m_{13}P^*_3 \right)<0, \quad 
J_{22}=-\frac{r_2}{k_2}P^*_2<0, \quad
J_{33}=-\frac{r_3}{k_3}P^*_3-\frac {m_{32}}{P_3^*}P_2^*<0,
$$
showing that it is always stable, when feasible.

In this case the two new equilibrium points arise
$I_2=(k_1,0,0)$, $I_3=(\alpha,0,\beta)$, with
$$
\alpha=\displaystyle\frac{k_1}{2}\left(1+\sqrt{1+\frac{4m_{13}k_3(r_3-m_{13})}{r_1r_3k_1}}\right), \quad
\beta=\displaystyle\frac{k_3}{r3}(r_3-m_{13}) ,
$$
feasible for
\begin{equation}\label{feas_I3}
r_3>m_{13}.
\end{equation}

At equilibrium $I_2$ we find the eigenvalues
$J_{11}=-r_1<0$, $J_{22}=r_2-m_{12}-m_{32}$, $J_{33}=r_3-m_{13}$.
Stability conditions are therefore
\begin{equation}\label{stab_I2}
r_2<m_{12}+m_{32}, \quad r_3<m_{13}.
\end{equation}

When $I_2$ is stable, the equilibria $I_3$ and coexistence are infeasible.
Thus at $r_2^{\dagger}=m_{12}+m_{32}$ and $r_3^{\dagger}=m_{13}$ there
are two transcritical bifurcations, the first one taking $I_2$ into coexistence,
the second one taking it into $I_3$.

At equilibrium $I_3$, the eigenvalues are
$$J_{11}=-\frac {r_1}{k_1} \alpha - m_{13} \frac {\beta}{\alpha }<0, \quad
J_{22}=r_2-m_{12}-m_{32}, \quad J_{33}=-\frac {r_3}{k_3}\beta<0.
$$
Stability is ensured by the first condition (\ref{stab_I2}).
Thus stability of $I_3$ also prevents feasibility of the coexistence equilibrium.
We have thus another transcritical bifurcation at $r_2^{\dagger}=m_{32}+m_{12}$
taking $I_3$ into coexistence.


\subsection{The model (\ref{n1})}

Here $m_{13}=m_{31}=m_{12}=m_{23}=0$.
This implies that at the origin one eigenvalue is $J_{33}=r_3>0$ so that this equilibrium
is unstable.

Coexistence can be stably attained, in view of the eigenvalues
$$
J_{11}=-\frac {r_1}{k_1}P_1^*<0, \quad
J_{22}=-\frac {r_2}{k_2}P_2^*-\frac {m_{21}}{P_2^*}P_1^*<0, \quad
J_{33}=-\frac {r_3}{k_3}P_3^*-\frac {m_{32}}{P_3^*}P_2^*<0, \quad
$$
at the levels
\begin{eqnarray*}
P_1^*= \frac {k_1}{r_1}\left( r_1-m_{21}\right), \quad
P_2^*=\frac {k_2}{2r_2}
\left[ r_2-m_{32}+ \sqrt{ \left( r_2-m_{32}\right)^2 +
\frac 4{k_2} r_2m_{21} P_1^*} \right] ,\\
P_3^*=\frac {k_3}{2r_3}\left[ r_3+ \sqrt{ r_3^2 +
\frac 4{k_3} r_3m_{32} P_2^*} \right] .
\end{eqnarray*}

In this case we find also the equilibrium $W_2=(0,0,k_3)$
and
$W_3=(0,P_2^+,P_3^+)$, with
$$
P_2^+=\frac {k_2}{r_2} (r_2-m_{32}), \quad
P_3^+=\frac {k_3}{2r_3}\left[ r_3 + \sqrt{ r_3^2 +
\frac 4{k_3} r_3m_{32} P_2^+} \right] .
$$

At $W_2$ the eigenvalues are $J_{11}=r_1-m_{21}$, $J_{22}=r_2-m_{32}$, $J_{33}=-r_3<0$,
providing stability when
\begin{equation}\label{stab_Q2}
r_1<m_{21}, \quad r_2<m_{32}.
\end{equation}

At $W_3$ the eigenvalues are $J_{11}=r_1-m_{21}$ and
$$J_{22}=-\frac {r_2}{k_2}P_2^+<0, \quad
J_{33}=-\frac {m_{32}}{P_3^+}P_2^+-\frac {r_3}{k_3}P_3^+<0,
$$
giving stability when
the first condition (\ref{stab_Q2}) holds.


\subsection{The model (\ref{n2})}
When $m_{13}=m_{31}=m_{12}=m_{32}=0$,
one eigenvalue at the origin is $J_{22}=r_2>0$, showing its instability.
Coexistence is allowed at the population values
\begin{eqnarray*}
P_1^*= \frac {k_1}{r_1}\left( r_1-m_{21}\right), \quad
P_3^*= \frac {k_3}{r_3}\left( r_3-m_{23}\right) ,\\
P_2^*=\frac {k_2}{2r_2}
\left[ r_2+ \sqrt{ r_2^2 +
\frac 4{k_2} r_2 (m_{21} P_1^*+m_{23} P_3^*)} \right] .
\end{eqnarray*}
This equilibrium is feasible for
\begin{equation}\label{feas_P*_n2}
r_1>m_{21}, \quad r_3>m_{23}.
\end{equation}
The eigenvalues are always negative, so that it is always stable, when feasible:
$$
J_{11}=-\frac {r_1}{k_1}P_1^*<0, \quad
J_{22}=-\frac {r_2}{k_2}P_2^*-\frac {m_{21}}{P_2^*}P_1^*-\frac {m_{23}}{P_2^*}P_3^*<0, \quad
J_{33}=-\frac {r_3}{k_3}P_3^*<0.
$$

We find also the equilibrium $X_1=(0,k_2,0)$, stable for
\begin{equation}\label{stab_X1}
r_2<m_{21}, \quad r_3<m_{23},
\end{equation}
and
more equilibria with either patch 1 or patch 3 empty, namely
$X_2=(P_1^X,P_2^X,0)$, with $P_1=P_1^*$,
$$
P_2^X=\frac {k_2}{2r_2}
\left[ r_2+ \sqrt{ r_2^2 + \frac 4{k_2} r_2 m_{21} P_1^X} \right] .
$$
and $Y_3=(0,P_2^Y,P_3^Y)$, $P_3^Y=P_3^*$,
$$
P_2^Y=\frac {k_2}{2r_2}
\left[ r_2+ \sqrt{ r_2^2 + \frac 4{k_2} r_2 m_{23} P_3^Y} \right].
$$

The eigenvalues at $X_2$ are
$$
J_{11}=-\frac {r_1}{k_1}P_1^X<0, \quad
J_{22}=-\frac {r_2}{k_2}P_2^X-\frac {m_{21}}{P_2^X}P_1^X<0, \quad
J_{33}=r_3-m_{23},
$$
giving stability for
\begin{equation}\label{feas_X2}
r_3<m_{23}.
\end{equation}

At $Y_3$ we have instead
$$
J_{11}=r_1-m_{21}, \quad
J_{22}=-\frac {r_2}{k_2}P_2^Y-\frac {m_{23}}{P_2^Y}P_3^Y<0, \quad
J_{33}=-\frac {r_3}{k_3}P_3^Y<0
$$
and consequently the stability conditions become
\begin{equation}\label{feas_Y3}
r_1<m_{21}.
\end{equation}


\subsection{The model (\ref{n3})}
Finally we consider $m_{13}=m_{31}=m_{21}=m_{23}=0$, with the eigenvalue $J_{11}=r_1>0$
at the origin, giving instability.

For coexistence we find the population levels
\begin{eqnarray*}
P_2^*= \frac {k_2}{r_2}\left( r_2-m_{12}-m_{32}\right), \\
P_1^*=\frac {k_1}{2r_1}
\left[ r_1+ \sqrt{ r_1^2 +
\frac 4{k_1} r_1 m_{12} P_2^*} \right] ,\\
P_3^*=\frac {k_3}{2r_3}
\left[ r_3+ \sqrt{ r_3^2 +
\frac 4{k_3} r_3 m_{32} P_3^*} \right] .
\end{eqnarray*}
It is feasible for
\begin{equation}\label{feas_P*_16}
r_2>m_{32}+m_{12}
\end{equation}
and when feasible it is always stable, since the eigenvalues are
$$
J_{11}=-\frac {r_1}{k_1}P_1^*-\frac {m_{12}}{P_1^*}P_2^*<0, \quad
J_{22}=-\frac {r_2}{k_2}P_2^*<0, \quad
J_{33}=-\frac {r_3}{k_3}P_3^*-\frac {m_{32}}{P_3^*}P_2^*<0.
$$

In addition to the origin and coexistence, we find $Z_1=(k_1,0,0)$, $Z_2=(0,0,k_3)$ and
$Z_3=(k_1,0,k_3)$. The former two are unstable, one eigenvalue is positive, respectively
$J_{11}=r_1>0$ and $J_{33}=r_3>0$. For the last equilibrium $Z_3$, the eigenvalues are
$J_{11}=-r_1<0$, $J_{22}=r_2-m_{12}-m_{32}$, $J_{11}=-r_3<0$ giving stability for
\begin{equation}\label{stab_Z3}
r_2< m_{12}+m_{32}.
\end{equation}


\newpage

\section{Discussion}

\subsection{The original model}
The findings of this paper show that coexistence can be attained always, for the general
model and for all the other models in which some interconnecting paths become unpracticable.
When it is feasible, and when the local stability analysis can be performed, it appears that
whenever feasible, the coexistence equilibrium is also locally asymptotically stable. We conjecture
that in such case it is also globally asymptotically stable, in view
of most of the other results, including the transcritical bifurcations found in some of the reduced
models.

Another good result from the conservationist point of view is that the ecosystem never disappears,
in almost all the cases in which
the stability conditions for the origin can be evaluated explicitly. In view of the logistic
growth assumption for the populations in each patch, we conjecture that this result holds true also
for the more interconnected models. A notable exception is given however by models (\ref{7}) and (\ref{7bis}),
see conditions (\ref{stab_orig_mod7bis}).

\subsection{The broken paths models}
For the models (\ref{2}), (\ref{3}), (\ref{1})
our results are again good from both the conservationist point of view as well as for the
development of human artifacts,
because it shows that in these cases some of the migration paths can be removed without harming too
much the whole ecosystem behavior. Its equilibria indeed remain the same of the original model (\ref{sist4}),
namely the origin and coexistence.
However some changes occur in the stability conditions
of these equilibria. Therefore changes leading to the situations modeled by systems
(\ref{2}), (\ref{3}), (\ref{1}) should be treated with care. It is also interesting to note that
these results hold for models where either one, two or three arcs are removed, therefore they really
depend on the configuration of the system, rather than on the number of allowed connecting paths
between patches.

Note that the common characteristic of these models, which is not shared by all the other ones,
is that there is always the possibility of cycling between all the patches, i.e. starting from
patch 1, say, to go to patch 2 and then 3 and finally returning to patch 1. Thus, when this cycle
can be performed, no equilibria other than survival in all patches is allowed, except possibly ecosystem
disappearance.

In general, in all other models the allowed system stable configurations, apart from origin and
coexistence, are those of the patches
with incoming paths. For instance in model (\ref{n3}) the population cannot survive only in patch 2, since it
contains only outgoing paths.

Models (\ref{7}) and (\ref{8}) are also interesting. The two interconnected patches are a possible system configuration
in the former but not in the latter. The reason, once more, is the fact that in model (\ref{7}) the two patches
have both incoming paths from the remaining patch, while in (\ref{8}) they are sources for the paths leading
to the remaining patch. This feature is also present in model (\ref{2bis}), where the interconnected patches cannot
be stably present in the equilibrium configurations because one of the interconnected patches is the origin
of an outgoing path. Therefore just the presence of just
one such outgoing path is enough to destabilize an interconnection. Conversely, in model (\ref{7bis}) the interconnected
patches are stable, because the remaining patch is a source of an outgoing flow.

The remaining models share another interesting property. In model (\ref{6}) there are two additional possible stable
configurations. The patch with both incoming paths, which is expected in view of the considerations holding for
the previous situations, and the configuration with this patch together with the intermediate node. In other
words only the patch from which both paths are outgoing cannot be present in the stable configurations. A
similar result holds in model (\ref{n1}), where only the source patch cannot be stable, with either the end sink
patch or both the intermediate and the sink patches give rise to stable configurations. Model (\ref{n2}) combines
these results, as the sink or either one of the other two patches can be stable. Model (\ref{n3}) again excludes the only
patch that is a source for both paths connecting it to the other two patches, while both these two patches
can be stable at the same time.

\end{document}